\documentclass{article}

\usepackage{arxiv}

\usepackage[utf8]{inputenc} % allow utf-8 input
\usepackage[T1]{fontenc}    % use 8-bit T1 fonts
\usepackage{hyperref}       % hyperlinks
\usepackage{url}            % simple URL typesetting
\usepackage{booktabs}       % professional-quality tables
\usepackage{amsfonts}       % blackboard math symbols
\usepackage{nicefrac}       % compact symbols for 1/2, etc.
\usepackage{microtype}      % microtypography
\usepackage{lipsum}

\usepackage[final]{graphicx}

\title{GENERALIZED COEFFICIENTS OF THE DIRICHLET SERIES}

\author{
 Kirill V.~Kapitonets \\
 BAUMAN MSTU, GRADUATE 1990\\
  MCC EuroChem\\
  Moscow\\
  Russian Federation\\
  \texttt{kkapitonets@live.com} \\
}

\begin{document}
\maketitle
\begin{abstract}
The paper considers a method for converting a divergent Dirichlet series into a convergent Dirichlet series by directly converting the coefficients of the original series $1\rightarrow\delta_{n}(s)$ for the Riemann Zeta function.
\par
In the first part of the paper, we study the properties of the coefficients ${\delta}^*_n$ of a finite Dirichlet series for approximating the Riemann Zeta function on the interval $\Delta{H}$.
\par
In general, the coefficients ${\delta}^*_n$ of a finite Dirichlet series are complex numbers.
\par
The dependence of the coefficients ${\delta}^*_n$ of a finite Dirichlet series on the ordinal number of the coefficient $n$ is established, which can be set by a sigmoid, and for each $N$ there is a single sigmoid $\hat{\delta}_n$ and a single interval $\Delta{H}$ for which the condition is satisfied
$$\Big|\sum\limits_{n}^{N}\{{\delta}^*_n- \hat\delta_n\}\Big| < \epsilon;$$
\par
The second part of the paper presents the results of using the sigmoid to calculate the values of the generalized coefficients $\delta_{n}(s)$ of the Dirichlet series for the Riemann Zeta function.
\par
For the accuracy of the calculation, $\log_{10}(1/\epsilon)$ values of the Riemann Zeta Function by summing the resulting convergent Dirichlet series, the power of the imaginary part $t=Im(s)$ is established.
\par
Presumably, the sigmoid can be used to asymptotically calculate the values of the analytical continuation of any function defined by the Dirichlet series.
\par
Presumably, for any divergent series to which the generalized summation method is applicable, it is possible to find a direct transformation of the coefficients of the divergent series, so that the resulting series with the transformed coefficients will converge to the same function as the series of transformed partial sums.
\end{abstract}
% keywords can be removed
\keywords{Dirichlet series, finite Dirichlet series, Riemann zeta function, generalized summation,  analytic continuation
}

%begin of text
\section{Introduction}
It is known \cite{BO} that the Riemann Zeta function is a function of a complex variable that is defined on the entire complex plane with a simple pole at $s = 1$ with residue $1$ and is an analytical continuation of the function of a complex variable\footnote{We deliberately changed the classical definition of the Riemann Zeta function given in the definition of the Riemann hypothesis, because we are not inclined to identify the whole function, which is the Riemann Zeta function, and the function that is defined in a bounded domain.}, which is defined in the right half-plane $Re (s) > 1$ by the Dirichlet series
\begin{equation}\label{dirichlet_zeta}\sum\limits_{n=1}^{\infty}\frac{1}{n^s}; \end{equation}
\par
where this series absolutely converges.
\par
In addition, the Riemann Zeta function satisfies the functional equation
\begin{equation}\label{zeta_func_eq}\pi^{-s/2}\Gamma(\frac{s}{2})\zeta(s)=\pi^{-(1-s)/2}\Gamma(\frac{1-s}{2})\zeta(1-s);\end{equation}
\par
Consequently, the values of the Riemann Zeta function can be easily determined on the entire complex plane through the values of the Dirichlet series in the right half-plane, with the exception of the critical strip  $0<\sigma<1$\footnote{Here and further we use the traditional notation of the complex variable in the theory of the Riemann Zeta function $s=\sigma+it$.}, where the Dirichlet series diverges for both the right and left sides of the functional equation.
\par
Thus, to determine the values of the Riemann Zeta function in the critical band, it is necessary to use special methods of generalized summation.
\par
Traditionally, in the theory of the Riemann Zeta function  \cite{TI}, the Euler-Maclaurin formula is used for the generalized summation of the Dirichlet series
\begin{equation}\label{eiler_macloren}\sum_{a<x}^b \phi(x)=\int_{a}^{b} \phi(x)dx+\int_{a}^{b}(x-[x]-\frac{1}{2}) \phi'(x)dx+(a-[a]-\frac{1}{2}) \phi(a)-(b-[b]-\frac{1}{2}) \phi(b);\end{equation}
\par
on the basis of which the integral definition of the Riemann Zeta function is obtained
\begin{equation}\label{zeta_int}\zeta(s) = s\int_{0}^{\infty}\frac{[x]-x}{x^{s+1}}dx;  0<\sigma<1;\end{equation}
\par
as well as a number of asymptotic formulas \cite{HA2}, including for determining the values of the Riemann Zeta function in the critical strip
\begin{equation}\label{zeta_asy_1}\zeta (s)=\lim \limits_{N\rightarrow +\infty}\left(\sum \limits_{n=1}^{N}{\frac {1}{n^{s}}}-{\frac {N^{1-s}}{1-s}}\right);\sigma>0;\end{equation}
\begin{equation}\label{zeta_asy_2}\zeta (s)=\lim \limits _{N\rightarrow +\infty }\left(\sum \limits _{n=1}^{N}{\frac {1}{n^{s}}}-{\frac {N^{1-s}}{1-s}}-{\frac {1}{2}}N^{-s}\right); \sigma>-1;\end{equation}
\begin{equation}\label{zeta_asy_3}\zeta (s)=\lim \limits _{N\rightarrow +\infty }\left(\sum \limits _{n=1}^{N}{\frac {1}{n^{s}}}-{\frac {N^{1-s}}{1-s}}-{\frac {1}{2}}N^{-s}+{\frac {1}{12}}sN^{-1-s}\right); \sigma>-2;\end{equation}
\par
Generalized summation methods \cite{HA2} operate with sequences of partial sums, for example, for the Cesaro method $(C, 1)$ for a divergent series
\begin{equation}\label{c1}\sum_{k=1}^{\infty} a_k;\end{equation}
\par
the sequence of partial sums of this series is considered
\begin{equation}\label{c2}s_m=\sum_{k=1}^{m} a_k;\end{equation}
\par
from which a sequence of average values is constructed
\begin{equation}\label{c3}t_m=\frac{s_0+s_1+s_2+...+s_m}{m+1};\end{equation}
\par
then the limit of these average values is considered
\begin{equation}\label{c4}t_m=\sum_{n=0}^{\infty} c_{mn}s_n\rightarrow s;\end{equation}
\par
It would be interesting to get a convergent series directly by converting $1\rightarrow\delta_{n}(s)$ coefficients of the original divergent series, for the Riemann Zeta function, such a series is the Dirichlet series.
\begin{equation}\label{zeta_delta_s}\zeta(s)=\sum_{n=1}^{\infty} \delta_{n}(s)n^{-s};\end{equation}
\par
Yuri Matiyasevich \cite{MA} set the task to determine the behavior of the coefficients of a finite Dirichlet series based on a sequence of non-trivial zeros of the Riemann Zeta function.
\begin{equation}\label{finite_dirihlet}\zeta(s)=\sum_{n=1}^{N} {\delta}^*_n n^{-s};\end{equation}
\par
As a result of the use of zeros of the Riemann Zeta function, the problem is reduced to the disclosure of the determinant by the column of terms of the Dirichlet series. But it is not the method of solving the problem that is important, but the result.
\par
As is known, the zeros of the Riemann Zeta function partially coincide with the zeros of the function
\begin{equation}\label{eta}\eta(s)=(1-2^{1-s})\zeta(s);\end{equation}
\par
which is set by the alternating Dirichlet series
\begin{equation}\label{dirichlet_eta}\sum\limits_{n=1}^{\infty}(-1)^{n-1}\frac{1}{n^s};\end{equation}
\par
which conventionally converges in the critical strip.
\par
Therefore, instead of the coefficients defining the Riemann Zeta function, the coefficients of the alternating Dirichlet series were obtained.
\begin{figure}[ht!]
\centering
\includegraphics[scale=0.5]{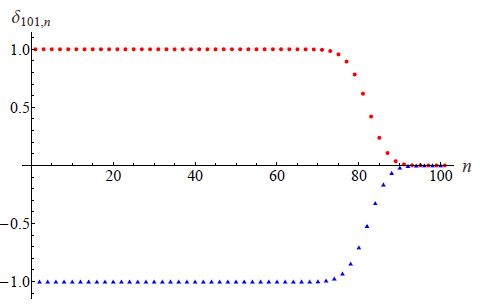}
\caption{Coefficients of a finite Dirichlet series \cite{MA}}
\label{fig:eta_cf}
\end{figure}
\par
Thus, the problem of obtaining the coefficients of a finite Dirichlet series for the Riemann Zeta function was not solved.
\par
Therefore, in our study, we set the task to calculate the coefficients ${\delta}^*_n$ of the finite Dirichlet series for the Riemann Zeta function (\ref{finite_dirihlet}) and in addition to find an expression for calculating the generalized coefficients $\delta_n (s)$ of the Dirichlet series for the Riemann Zeta function (\ref{zeta_delta_s}).
\par
\section{Calculation of the coefficients of a finite Dirichlet series}
To calculate the coefficients of a finite Dirichlet series for the Riemann Zeta function (\ref{finite_dirihlet}), we solve a system of linear equations (\ref{linear_eq}).
\begin{equation}\label{linear_eq}Ax=b;\end{equation}
\par
As coefficients $\{a_{mn}\}$ , we will use the values of the corresponding terms of the Dirichlet series (\ref{dirichlet_e}) calculated with a given accuracy $log_{10}(1/\epsilon)$\footnote{For the calculation, we use the Python mathematical package npmath.} (\ref{e1}).
\begin{equation}\label{dirichlet_e} a_{mn}=\hat n_m^{-s}; \end{equation}
\begin{equation}\label{e1}|\hat n_m^{-s}-n_m^{-s}|<\epsilon; \end{equation}
\par
In order for the system of equations to be non-degenerate, we will use as the vector of free terms $\{b_m\}$ of the values of the Riemann Zeta function in the critical strip other than zero (\ref{zeta_e}) calculated with a given accuracy $log_{10}(1/\epsilon)$ (\ref{e2}). 
\begin{equation}\label{zeta_e}b_m=\hat\zeta(s_m); \end{equation}
\begin{equation}\label{e2}|\hat\zeta(s_m)-\zeta(s_m)|<\epsilon; \end{equation}
\par
Assumption 1.
Generalized Dirichlet series coefficients for the Riemann Zeta function $\delta_n (s)$ are real numbers and depend on the value of the complex variable of the Riemann Zeta function.
\begin{equation}\label{delta_s_r}\delta_n(s) \subset \mathbb{R};\end{equation}
\par
Statement 1.
The coefficients of a finite Dirichlet series ${\delta}^*_n$ are complex numbers
\begin{equation}\label{delta_s_c}{\delta}^*_n \subset \mathbb{C};\end{equation}
\par
Obviously, in accordance with assumption 1, it is possible to find real coefficients for a finite Dirichlet series $\delta_{mn}$ corresponding to individual values of $s_m$.
\begin{equation}\label{finite_dirihlet_re}\hat \zeta(s_m)=\sum_{n=1}^{N} \delta_{mn} \hat n_m^{-s};\end{equation}
\begin{equation}\label{delta_mn_re}\delta_{mn} \subset \mathbb{R};\end{equation}
\par
On the other hand, the same values are $\hat\zeta(s_m)$ can be written in terms of the coefficients of a finite Dirichlet series ${\delta}^*_n$, which are the solution of a system of linear equations (\ref{linear_eq})
\begin{equation}\label{finite_dirihlet_im}\hat \zeta(s_m)=\sum_{n=1}^{N} {\delta}^*_n \hat n_m^{-s};\end{equation}
\par
Obviously, the sum (\ref{zeta_sum1}) of all the values of the Riemann Zeta function will not depend on which coefficients we expressed each individual value through
\begin{equation}\label{zeta_sum1}\sum\limits_{m}^{N} \hat \zeta(s_m);\end{equation}
\par
Therefore
\begin{equation}\label{zeta_sum2}\sum\limits_{m}^{N} \sum\limits_{n}^{N} {\delta}^*_n n^{-s}_m = \sum\limits_{m}^{N} \sum\limits_{n}^{N} \delta_{mn} n^{-s}_m;\end{equation}
\par
Let's change the order of summation and extract the coefficients ${\delta}^*_n$ of the finite Dirichlet series from under the sign of the sum, since from the solution of the equation system we will get the same coefficients for all given $s_m$.
\begin{equation}\label{cf2}\sum\limits_{n}^{N} {\delta}^*_n \sum\limits_{m}^{N}  n^{-s}_m = \sum\limits_{n}^{N} \sum\limits_{m}^{N} \delta_{mn} n^{-s}_m;\end{equation}
\par
Comparing the expressions under the sign of the first sum, we get
\begin{equation}\label{cf3}{\delta}^*_n = \frac{\sum\limits_{m}^{N} \delta_{mn} n^{-s}_m}{\sum\limits_{m}^{N}n^{-s}_m};\end{equation}
\par
Obviously, in the general case, the values of ${\delta}^*_n$ are complex numbers.
\par
Corollary 1.
The coefficients of the finite Dirichlet series ${\delta}^*_n$ are the average values of the real coefficients of the finite Dirichlet series $\delta_{mn}$ for the given values $s_m$.
\par
Statement 2.
The maximum value of $t_m=Im(s_m)$ is limited by the number of coefficients of a finite Dirichlet series
\begin{equation}\label{big_n}max(\frac{t_m}{\pi}) < N;\end{equation}
\par
The asymptotic expression \cite{HA1} obtained based on the Euler-Maclaurin formula
\begin{equation}\label{zeta_eq_1}\zeta(s)=\sum_{n\le x}{\frac{1}{n^s}}-\frac{x^{1-s}}{1-s}+\mathcal{O}(x^{-\sigma}); 
\sigma>0; |t|< 2\pi x/C; C > 1\end{equation}
\par
determines that to calculate the value of the Riemann Zeta function, it is necessary to take at least $\lfloor t/\pi \rfloor$ terms of the Dirichlet series.
\par
Therefore, to calculate the value of the Riemann Zeta function, it is necessary to know at least $\lfloor t/\pi \rfloor$ of the generalized Dirichlet series coefficients.
\par
Corollary 2.
The coefficients of the finite Dirichlet series ${\delta}^*_n$, depend simultaneously on the number of coefficients $N$ and the range of values $t_m=Im (s_m)$
\par
\section{Calculation results}
Analysis of the calculation of ${\delta}^*_n$ for a given value of $Re(s_m)$ and linear change $t_m=Im (s_m)$
\begin{equation}\label{im_s_m}t_m = t_1 + (m-1) \Delta t; m = 1..N;\end{equation}
\par
showed that the result depends both on the selected accuracy of calculating $log_{10}(1/\epsilon)$, and on the selected values of $t_1$ and $\Delta t$.
\par
\begin{figure}[ht!]
\centering
\includegraphics[scale=0.5]{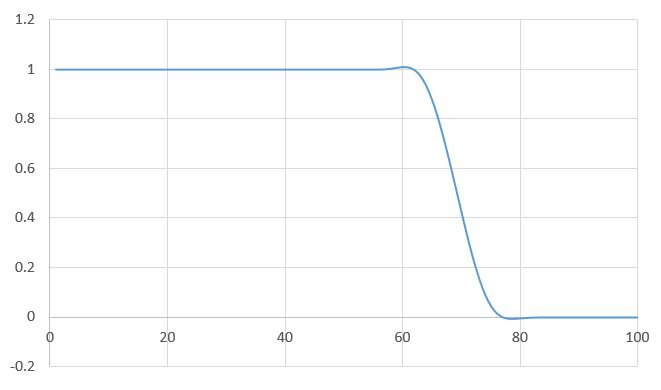}
\caption{The real part of ${\delta}^*_n$: $t_1=188.4955592$, $\Delta{t}=0.628318531$, $log_{10}(1/\epsilon)=100$}
\label{fig:re_100}
\end{figure}
\begin{figure}[ht!]
\centering
\includegraphics[scale=0.5]{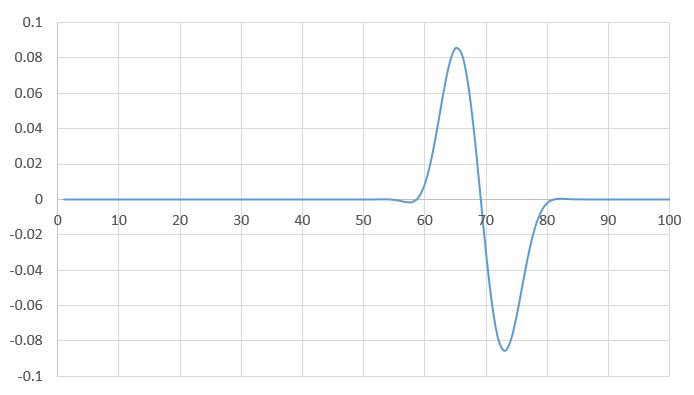}
\caption{Imaginary part of ${\delta}^*_n$}
\label{fig:im_100}
\end{figure}

Result 1.
The stable dependence of the coefficients of the finite Dirichlet series ${\delta}^*_n$ on the coefficient number $n$ occurs at certain ratios of the number of coefficients $N$, the initial $t_1$ and the final $t_1+(N-1) \Delta t$ values $t_m=Im(s_m)$ Fig. \ref{fig:re_100}, \ref{fig:im_100}.
\par
Any slight change in the specified parameters leads to a violation of the dependence of the value of ${\delta}^*_n$ on the coefficient number $n$.
\begin{figure}[ht!]
\centering
\includegraphics[scale=0.5]{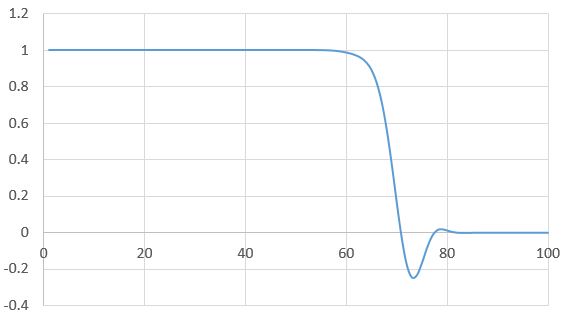}
\caption{$Re({\delta}^*_n)$, left offset: $N=100$, $t_1=157.0796327$, $\Delta{t}=0.785398163$, $log_{10}(1/\epsilon)=100$}
\label{fig:re_100_4}
\end{figure}
\begin{figure}[ht!]
\centering
\includegraphics[scale=0.5]{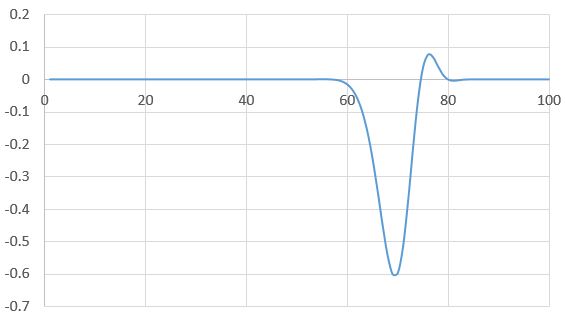}
\caption{$Im({\delta}^*_n)$, left offset}
\label{fig:im_100_4}
\end{figure}
\par
If we reduce the minimum value of $t_1=Im (s_1)$, then the values of ${\delta}^*_n$ are deformed to the left Fig. \ref{fig:re_100_4}, \ref{fig:im_100_4}.
\begin{figure}[ht!]
\centering
\includegraphics[scale=0.5]{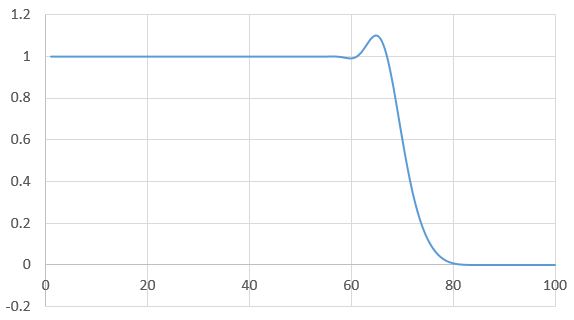}
\caption{$Re({\delta}^*_n)$, right offset: $t_1=209.4395102$, $\Delta{t}=0.523598776$, $log_{10}(1/\epsilon)=100$}
\label{fig:re_100_6}
\end{figure}
\begin{figure}[ht!]
\centering
\includegraphics[scale=0.5]{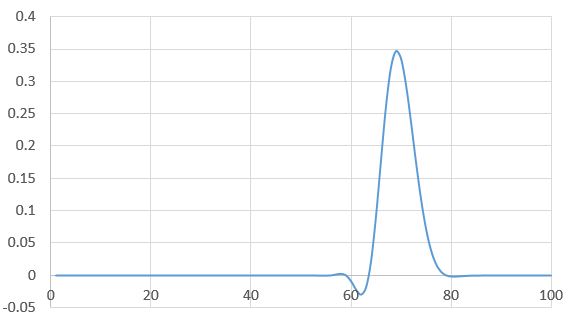}
\caption{$Im({\delta}^*_n)$, offset to the right}
\label{fig:im_100_6}
\end{figure}
\par
Accordingly, if we increase the minimum value of $t_1=Im (s_1)$, then the values of ${\delta}^*_n$ are deformed to the right Fig. \ref{fig:re_100_6}, \ref{fig:im_100_6}.
\begin{figure}[ht!]
\centering
\includegraphics[scale=0.5]{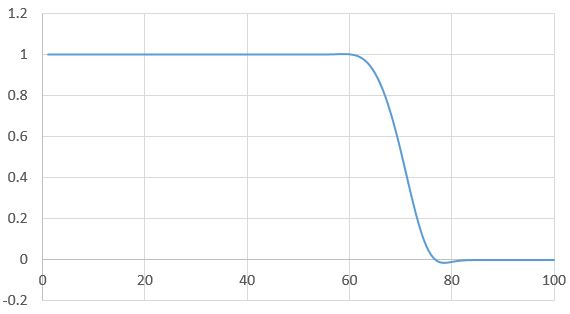}
\caption{$Re({\delta}^*_n)$, a slight decrease in accuracy: $t_1=188.4955592$, $\Delta{t}=0.628318531$, $log_{10}(1/\epsilon)=90$}
\label{fig:re_100_10}
\end{figure}
\begin{figure}[ht!]
\centering
\includegraphics[scale=0.5]{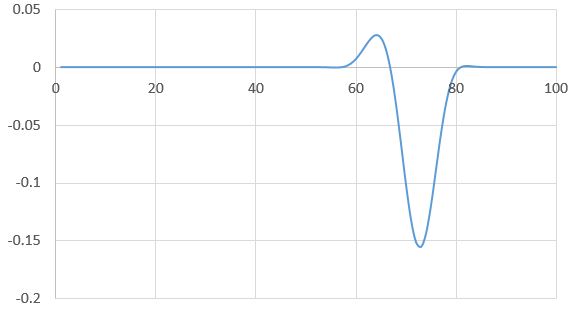}
\caption{$Im({\delta}^*_n)$, a slight decrease in accuracy}
\label{fig:im_100_10}
\end{figure}
\par
\begin{figure}[ht!]
\centering
\includegraphics[scale=0.5]{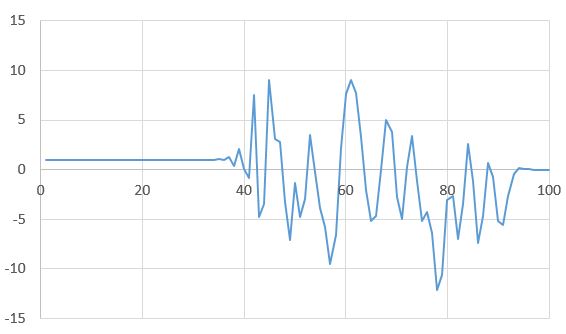}
\caption{$Re({\delta}^*_n)$, a significant decrease in accuracy: $t_1=188.4955592$, $\Delta{t}=0.628318531$, $log_{10}(1/\epsilon)=50$}
\label{fig:re_100_8}
\end{figure}
\begin{figure}[ht!]
\centering
\includegraphics[scale=0.5]{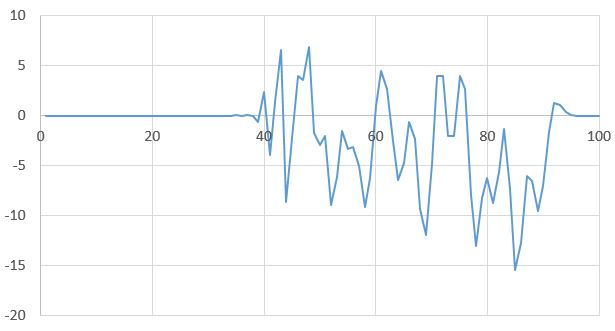}
\caption{$Im({\delta}^*_n)$, a significant decrease in accuracy}
\label{fig:im_100_8}
\end{figure}

If we slightly reduce the accuracy of the calculations $log_{10}(1/\epsilon)$ (\ref{e1}, \ref{e2}), then the values of ${\delta}^*_n$ are first deformed Fig. \ref{fig:re_100_10}, \ref{fig:im_100_10}, and with a significant decrease in the accuracy of calculations, the dependence of ${\delta}^*_n$ on the number disappears altogether Fig. \ref{fig:re_100_8}, \ref{fig:im_100_8}.
\par
Therefore, the condition for the stable dependence of the coefficients of the finite Dirichlet series ${\delta}^*_n$ on the coefficient number $n$ can be formulated as follows
\begin{equation}\label{delta_n_im_sum}\Big|\sum\limits_{n}^{N} Im({\delta}^*_n)\Big| < \epsilon;\end{equation}
\par
\begin{figure}[ht!]
\centering
\includegraphics[scale=0.5]{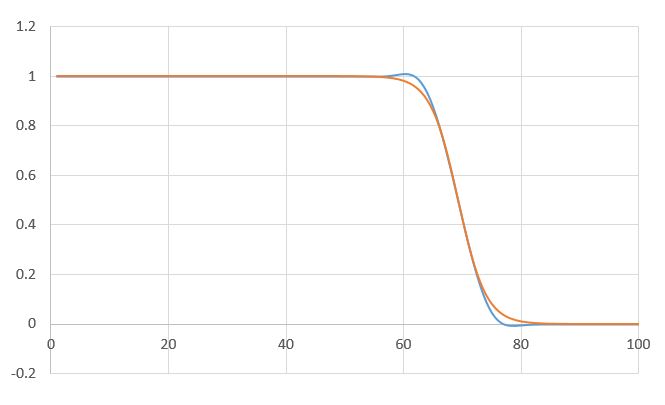}
\caption{Sigmoid and real values of $\hat\delta_n$}
\label{fig:re_100_re}
\end{figure}

Result 2.
For each stable dependence, Fig. \ref{fig:re_100} coefficients of the finite Dirichlet series ${\delta}^*_n$ from the coefficient number $n$ there is a sigmoid (\ref{hat_delta_n}) that defines a set of real values $\hat \delta_n$ that correspond to the values of the real part of the coefficients of the finite Dirichlet series Fig. \ref{fig:re_100_re}.
\begin{equation}\label{hat_delta_n}\hat\delta_n= \frac{1}{1+\exp\Big(\frac{n-A}{B}\Big)};\end{equation}
\par
Result 3.
The coefficient of the sigmoid $A$ corresponds to the average value of $\hat{n}^*$ of the coefficient number ${\delta}^*_n$, at which its value would be equal to 1/2.
\par
This value is easily obtained by linear interpolation
\par
let
\begin{equation}\label{a}{\delta}^*_a >1/2; {\delta}^*_b <1/2; b=a+1;\end{equation}
\par
then
\begin{equation}\label{hat_n_f}\hat{n}^*=a+\frac{1/2-{\delta}^*_b}{{\delta}^*_a-{\delta}^*_b}\end{equation}
\par
Result 4.
The scale factor of the sigmoid $B$ depends on the average value of $\hat{n}^*$ of the coefficient number ${\delta}^*_n$, at which its value would be equal to 1/2 and from the number of coefficients $N$ and does not depend on the range of values $\Delta{H}$.
\begin{equation}\label{b}B^2=\hat{n}^*-\frac{2N}{\pi};\end{equation}
\par
At first glance, it seems that to fulfill the condition (\ref{delta_n_im_sum}) of a stable dependence of the coefficients of a finite Dirichlet series ${\delta}^*_n$ on the coefficient number $n$, it is enough to choose such an average value $\hat{t}_m=Im (\hat{s}_m)$
\begin{equation}\label{hat_t}\hat{t}_m = t_1+\frac{(N-1)\Delta{t}}{2};\end{equation}
\par
to meet the condition
\begin{equation}\label{hat_n_f_eq_hat_n}\hat{n}^*=\hat{n};\end{equation}
\par
where
\begin{equation}\label{hat_n}\hat{n}=\frac{\hat{t}_m}{\pi} = \frac{t_1}{\pi}+\frac{(N-1)\Delta{t}}{2\pi};\end{equation}
\par
\begin{figure}[ht!]
\centering
\includegraphics[scale=0.5]{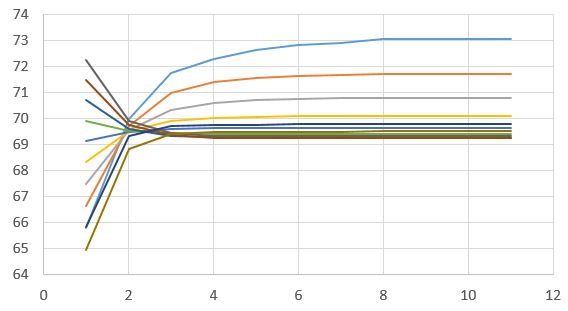}
\caption{Dependence of $\hat{n}^*$ on the initial value of $t_1$}
\label{fig:delta_h_m}
\end{figure}

Nevertheless, it turned out that the condition (\ref{hat_n_f_eq_hat_n}) can be fulfilled for any $t_1$ for which there is a stable dependence of the coefficients of the finite Dirichlet series ${\delta}^*_n$ on the coefficient number $n$, but at the same time $\hat{n}^*$ will converge to different values Fig. \ref{fig:delta_h_m}.
\par
\begin{figure}[ht!]
\centering
\includegraphics[scale=0.5]{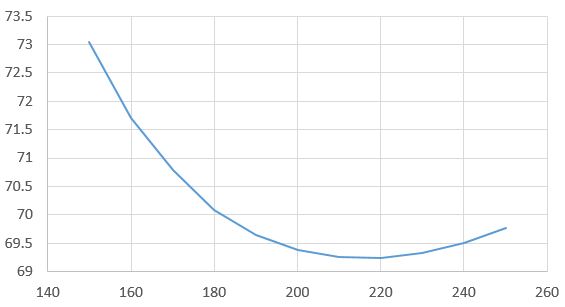}
\caption{Dependence of the average value of $\hat{n}$ on $t_1$}
\label{fig:hat_n_m}
\end{figure}

Result 5.
There is an average ratio of $\hat{t}_m/t_1$ Fig. \ref{fig:hat_n_m} when the condition is met
\begin{equation}\label{delta_n_sum}\Big|\sum\limits_{n}^{N}\{{\delta}^*_n- \hat\delta_n\}\Big| < \epsilon;\end{equation}
Corollary 3.
For each $N$, there is a minimum value of $\hat{n}^*$, at which the condition (\ref{delta_n_im_sum}) is met.
\par
Corollary 4.
For each $N$, there is a single sigmoid (\ref{hat_delta_n}) for which the condition (\ref{delta_n_sum}) is met.
\par
\section{Calculation of generalized coefficients of the Dirichlet series}
It is obvious that the sigmoid (\ref{hat_delta_n}) is a suitable expression for calculating the generalized coefficients $\delta_n (s)$ of the Dirichlet series (\ref{zeta_delta_s}) for the Riemann Zeta function.

When moving from a finite Dirichlet series (\ref{finite_dirihlet}) to an infinite one (\ref{zeta_delta_s}), it is logical to use the value $t=Im(s)$  to calculate the average number of $\hat{n}$
\par
\begin{equation}\label{im_s}\hat{n}=\frac{t}{\pi};\end{equation}
\par
Then the generalized coefficients $\delta_n(s)$ of the Dirichlet series for the Riemann Zeta function can be calculated as follows
\begin{equation}\label{zeta_delta_t}\delta_n(s)=\frac{1}{1+exp(\frac{n-t/\pi}{B(s)})};\end{equation}
\par
Obviously, we can no longer use the expression (\ref{b}) for the scale factor $B(s)$.
\par
Assumption 2.
There is an average value of the scale factor $\hat{B}(s)$ at which the condition is met
\begin{equation}\label{e3}|\zeta(s)-\sum_{n=1}^{N} \delta_{n}(s)n^{-s}|<\epsilon; \end{equation}
\begin{figure}[ht!]
\centering
\includegraphics[scale=0.5]{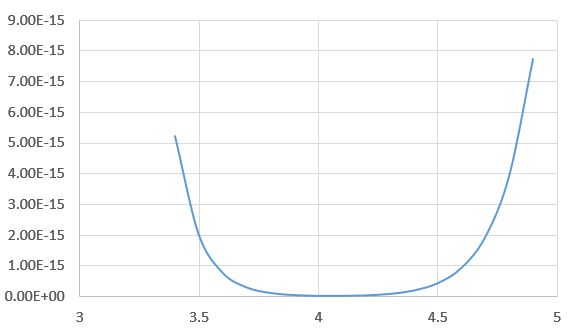}
\caption{Dependence of $\epsilon$ on $B$: s = $0.5+1000i$; $\hat{B}=4.05968$}
\label{fig:epsilon_1000}
\end{figure}
\par
Indeed, if the sigmoid (\ref{zeta_delta_t}) is very steep, then $n^{- s}$ will decrease too quickly and the sum (\ref{zeta_delta_s}) will not have time to gain the necessary value, on the other hand, if the sigmoid is too flat, then $n^{-s}$ will begin to decrease too early and the sum will not gain the necessary value again, although there may be some compensating effect of the slow decrease of $n^{-s}$.
\par
Nevertheless, the calculations show that the average value of $\hat{B}$ exists Fig. \ref{fig:epsilon_1000} and it depends on both $\sigma=Re (s)$ and $t=Im (s)$.
\par
\begin{figure}[ht!]
\centering
\includegraphics[scale=0.5]{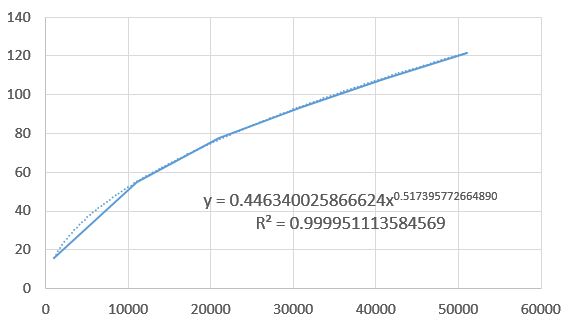}
\caption{Dependency $\log_{10}(1/\epsilon)$ from $t=Im (s)$}
\label{fig:epsilon_t}
\end{figure}

Result 6.
The values of the Riemann Zeta function obtained using generalized coefficients (\ref{zeta_delta_t}) can be calculated with an asymptotic\footnote{The accuracy of calculations (\ref{e3}) increases with the growth of the imaginary part of the complex variable $t=Im(s)$ Fig. \ref{fig:epsilon_t}} accuracy for any values of the complex variable Riemann Zeta function both in the right half-plane, where the Dirichlet series absolutely converges, and in the left half-plane, including in the critical band, where the Dirichlet series diverges, except for the real axis\footnote{It is not possible to calculate the coefficients $\delta_n (s)$ on the real axis.}
\par
\begin{figure}[ht!]
\centering
\includegraphics[scale=0.5]{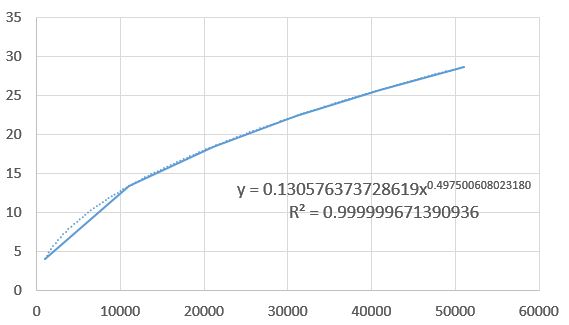}
\caption{Dependence of the proportionality coefficient $\hat{B}(s)$ from $t=Im (s)$: $\sigma=0.5$}
\label{fig:b_t}
\end{figure}

Result 7.
The average value of $\hat{B}$ is in power dependence (\ref{b_t}) on the imaginary part of the complex variable $t=Im(s)$ Fig. \ref{fig:b_t}
\begin{equation}\label{b_t}\hat{B}(\sigma+it)=C(\sigma)t^{D(\sigma)}; \end{equation}
\par
\begin{figure}[ht!]
\centering
\includegraphics[scale=0.5]{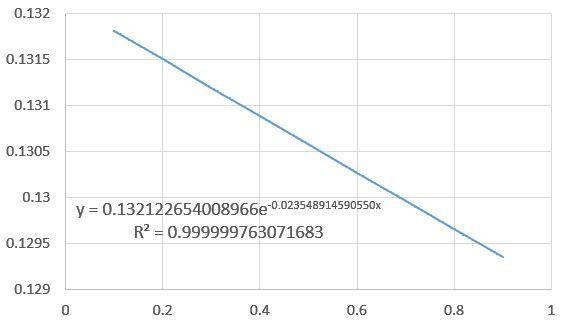}
\caption{The dependence of the coefficient C on $\sigma=Re (s)$}
\label{fig:c_sigma}
\end{figure}
\begin{figure}[ht!]
\centering
\includegraphics[scale=0.5]{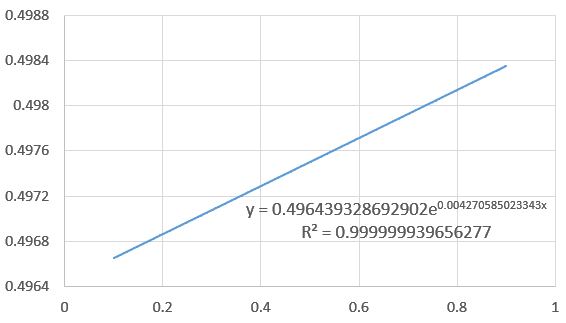}
\caption{The dependence of the coefficient D on $\sigma=Re (s)$}
\label{fig:d_sigma}
\end{figure}

Result 8.
The coefficients (\ref{b_t}) are in the inverse exponential \footnote{The dependence cannot be linear, because the proportionality coefficient of the sigmoid must be greater than zero.} dependencies on the real part of the complex variable $\sigma=Im(s)$ Fig. \ref{fig:c_sigma}, \ref{fig:d_sigma}.
\par
\begin{figure}[ht!]
\centering
\includegraphics[scale=0.5]{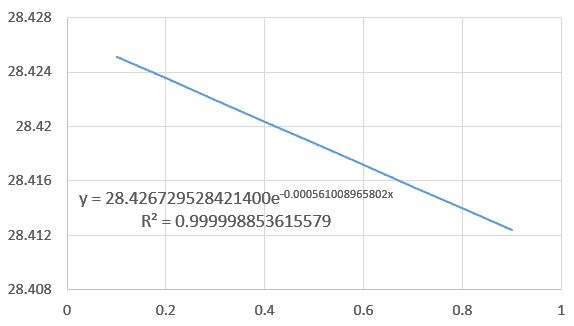}
\caption{Dependence of the proportionality coefficient $\hat{B}(s)$ from $\sigma=Re (s)$: $t=50000$}
\label{fig:b_sigma}
\end{figure}

Result 9.
The average value of $\hat{B}$ is exponentially dependent on the real part of the complex variable $\sigma=Im (s)$ Fig. \ref{fig:b_sigma}.
\par
\section{Conclusions}
On the one hand, there are many different methods for approximating functions, all of them are related to finding the coefficients of a linear combination of some functions\footnote{These functions must form an orthogonal basis or they must be linearly independent in some orthogonal basis}.
\par
On the other hand, there are many different methods of generalized summation of divergent series.
\par
But whatever method we use to find an infinite convergent series that defines the entire function of a complex variable, it will converge to the same values.
\par
In this paper, we tried to determine what the coefficients of the convergent Dirichlet series look like for the Riemann Zeta function and to determine a method for directly converting (\ref{1_to_delta}) the coefficients of the Dirichlet series, so that it becomes a convergent series from a divergent series.
\begin{equation}\label{1_to_delta}1\rightarrow\delta_{n}(s);\end{equation}
\par
It is obvious that we have obtained a convergent series, since the sigmoid (\ref{zeta_delta_t}) cuts off some of the terms of the infinite series, making their values negligible.
\par
On the other hand, we have shown that the condition (\ref{e3}) is asymptotically satisfied (Fig. \ref{fig:b_t}) for any values of the complex variable\footnote {Except for real values, for which it is not possible to construct a sigmoid.}.
\par
Therefore, we can assume that the resulting convergent series asymptotically, with the growth of the imaginary part of the complex variable $t=Im(s)$, converges to the values of the Riemann Zeta function.
\par
Assumption 3.
The $1\rightarrow\delta_{n}(s)$ transformation given by the sigmoid (\ref{zeta_delta_t}) transforms a Dirichlet series with coefficients 1, which diverges at $\sigma=Re (s)<1$ into a convergent series that asymptotically determines the values of the Riemann Zeta function.
\par
We can test our assumption by substituting the resulting convergent series into the functional equation (\ref{zeta_func_eq2})
\begin{equation}\label{zeta_func_eq2}\zeta(s)=\chi(s)\zeta(1-s); \end{equation}
\par
where
\begin{equation}\label{chi_eq}\chi(s)=\frac{(2\pi)^s}{2\Gamma(s)\cos(\large\frac{\pi s}{2})}=2^s\pi^{s-1}\sin(\frac{\pi s}{2})\Gamma(1-s)=\pi^{s-\frac{1}{2}}\frac{\Gamma(\frac{1-s}{2})}{\Gamma(\frac{s}{2})};\end{equation}
\par
Then obviously we must demand the fulfillment of the condition
\begin{equation}\label{e4}|\sum_{n=1}^{N} \delta_{n}(s)n^{-s}-\chi(s)\sum_{n=1}^{N} \delta_{n}(s)n^{s-1}| <\epsilon; \end{equation}
\par
Which we can rewrite as follows
\begin{equation}\label{e5}|\sum_{n=1}^{N} \delta_{n}(s)\{n^{-s}-\chi(s)n^{s-1}\}| <\epsilon; \end{equation}
\par
Let us first consider the divergent sum
\begin{equation}\label{e6}\sum_{n=1}^{N}\{n^{-s}-\chi(s)n^{s-1}\}; \end{equation}
\par
Partial sums of the specified series behave in no less strange way than partial sums of the Dirichlet series.
\par
At first, they form segments of a sufficiently large size, but then these segments become smaller and smaller, but the strangeness lies in something else.
\par
\begin{figure}[ht!]
\centering
\includegraphics[scale=0.5]{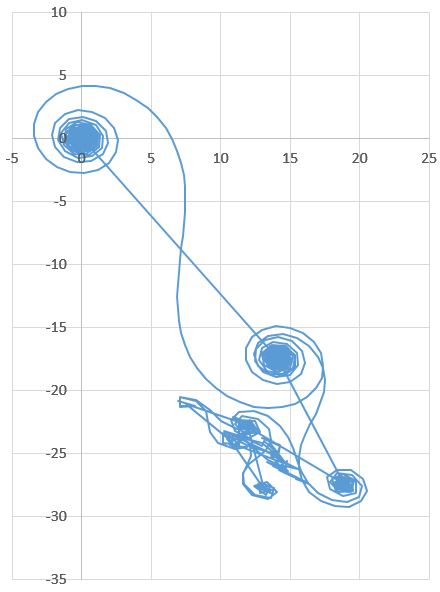}
\caption{Partial sums of a series (\ref{e6})}
\label{fig:spiral_1}
\end{figure}

Small segments are twisted around large segments and boomerang back to the origin of the complex plane Fig. \ref{fig:spiral_1}.
\par
\begin{figure}[ht!]
\centering
\includegraphics[scale=0.5]{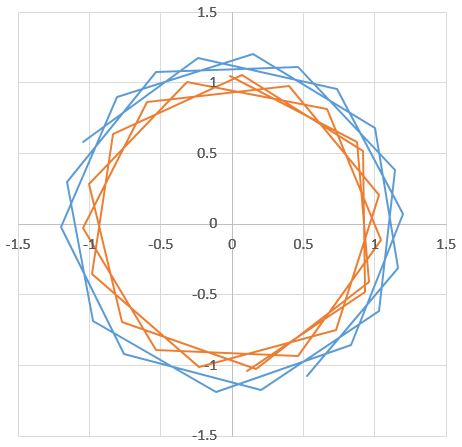}
\caption{Divergent spiral of partial sums of a series (\ref{e6})}
\label{fig:spiral_2}
\end{figure}

It is obvious that small segments are twisted around the origin of the coordinates of the complex plane Fig. \ref{fig:spiral_2} into a diverging spiral (the inner spiral is formed by segments with smaller numbers, and the outer spiral is formed by segments with larger numbers).
\par
\begin{figure}[ht!]
\centering
\includegraphics[scale=0.5]{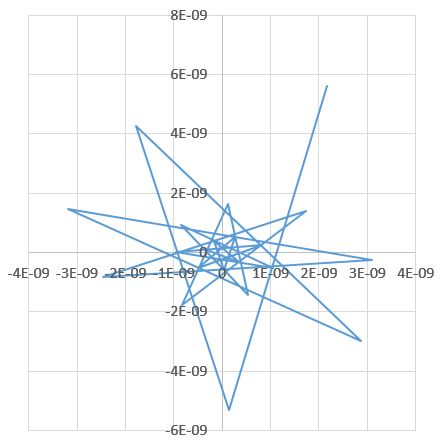}
\caption{Convergent spiral of partial sums of a series (\ref{e5})}
\label{fig:spiral_3}
\end{figure}

It is obvious that the sigmoid transforms a diverging spiral into a spiral converging to the origin of the coordinates of the complex plane Fig. \ref{fig:spiral_3} (the sigmoid affects segments that twist around the origin).
\par
Assumption 4.
The transformation $1\rightarrow\delta_{n} (s)$, given by the sigmoid (\ref{zeta_delta_t}), is a universal transformation of a divergent series of a function of a complex variable into a convergent series.
\par
The obtained results, of course, require a formal justification, but they show in what form it is necessary to look for functions of direct transformation of the coefficients of a divergent series in order to obtain a convergent series.
\par
But the most important result that we want to extract from the work done is the understanding that every whole function must be determined by a convergent series, if for some reason a divergent series is obtained, then it is undoubtedly necessary to look for the form of a convergent series, and even if generalized summation methods are used, then as such a convergent series it is necessary to consider the transformed series of partial sums, and not the original divergent series, attributing to it the sum obtained by generalized summation methods.
\par
Assumption 5.
For each divergent series, if the generalized summation method is applied to it, and this generalized summation method transforms the divergent series into a convergent series of transformed partial sums that defines an whole function, there is a method for directly converting the coefficients of the original divergent series, so that the series with the transformed coefficients will converge to the same integer function as the convergent series of transformed partial sums.

%end of text

\bibliographystyle{unsrt}  


\begin{thebibliography}{1}

\bibitem{MA}Beliakov, G. and Matiyasevich, Y. 2014, Approximation of Riemann's Zeta Function by Finite Dirichlet Series: A Multiprecision Numerical Approach. Experimental Mathematics. 24. 10.1080/10586458.2014.976801. 
\bibitem{HA1}Hardy, G. H. and J. E. Littlewood, 1921, The Zeros of Riemann's Zeta Function on the Critical Line. Mathematische Zeitschrift
\bibitem{HA2}G.H. Hardy, 1949, Divergent series, Oxford At The Clarendon Press, Available at:\\ https://archive.org/details/DivergentSeries
\bibitem{TI}Titchmarsh, E.C. (1988) The Theory of the Riemann Zeta Function. Oxford University Press, Oxford.
\bibitem{BO}Bombieri E., 2000, The Riemann Hypothesis - official problem description (PDF), Clay Mathematics Institute, retrieved 2008-10-25 Reprinted in (Borwein et al. 2008), Available at:\\ https://claymath.org/sites/default/files/official\_problem\_description.pdf

\end{thebibliography}

\end{document}